\documentclass[11pt]{article}
\usepackage{geometry}                % See geometry.pdf to learn the layout options. There are lots.
\usepackage{graphicx}
\usepackage{amssymb}
\usepackage{epstopdf}
\usepackage{amsmath,amsthm,amscd,amssymb}
\usepackage{latexsym}
\usepackage[colorlinks,citecolor=red,pagebackref,hypertexnames=false]{hyperref}

\numberwithin{equation}{section}

\theoremstyle{plain}
\newtheorem{theorem}{Theorem}[section]

\newtheorem{corollary}[theorem]{Corollary}
\newtheorem{proposition}[theorem]{Proposition}

\theoremstyle{definition}
\newtheorem{definition}[theorem]{Definition}

\theoremstyle{remark}

\newtheorem{case[theorem]}{Case}

\title{Direction sets, Lipschitz graphs and density} 
\author{Alex Iosevich and Jonathan Pakianathan}
\begin{document}
\maketitle

\begin{abstract} We consider the direction set determined by various subsets $E$ of Euclidean space and show that there is a trichotomy: Either  (i) The subset is the graph of a Lipschitz function and the direction set is not dense in the sphere, (ii) The subset is the graph of a non-Lipschitz function and the direction set is dense but not everything, or (iii) The subset is not a graph (in a suitable sense) and every direction is determined by the set.
We then explore a variety of results based on this trichotomy under additional assumptions on the set $E$.
\end{abstract}                

\section{Introduction} 

The purpose of this paper is study directions sets determined by subsets of the Euclidean space. Informally, direction sets consists of direction vectors determined by pairs of vectors from a given set. More precisely, we have the following definitions. 

\begin{definition} Fix a subset $E \subseteq \mathbb{R}^d, d \geq 2$. The (oriented) direction set $\tilde{D}(E)$ determined by $E$ is the set 
$$\tilde{D}(E)=\left\{ \frac{y-x}{|y-x|}: x, y \in E, x \neq y \right\} \subseteq S^{d-1}.$$ 

Note that this consists of the unit direction vector of the ray from $x$ to $y$ as $x, y$ range over pairs of distinct elements of $E$.  We say that $E$ determines all directions, or a dense subset of directions, respectively, if $\tilde{D}(E)=S^{d-1}$ or $\tilde{D}(E)$ is a dense subset of $S^{d-1}$, respectively.  Note also that when $|E| \leq 1$ then $\tilde{D}(E)=\emptyset$, but is nonempty otherwise.
\end{definition} 

It was shown by the first listed author, Mourgoglou and Senger (\cite{IMS12}; see also \cite{Mat95}, Theorem 10.11) that if the Hausdorff dimension of $E \subset {\Bbb R}^d$, $d \ge 2$, is greater than $d-1$, then the $(d-1)$-dimensional Lebesgue measure of $\tilde{D}(E)$, viewed as a subset of $S^{d-1}$, is positive. They also obtained a rather precise description of the distribution of these directions. In this paper we turn in a slightly different direction and obtain rather comprehensive qualitative information about the structure of subsets of ${\Bbb R}^d$ for which the direction set is not dense in the sphere. 

It is often convenient to consider unoriented directions, which are elements of $S^{d-1}$ taken modulo the antipodal action $a(x)=-x$. It is well known that the quotient of the sphere $S^{d-1}$ when antipodal points are identified is the projective space $\mathbb{R}P^{d-1}$ whose elements can be thought of as either pairs of antipodal points on the sphere, $\pm \hat{u}$ or lines through the origin in $\mathbb{R}^d$. The natural map $\pi: S^{d-1} \to RP^{d-1}$ is a double cover map such that for each line $L \in RP^{d-1}$, $\pi^{-1}(L)$ is the pair of antipodal points where the line $L$ intersects the sphere.

\begin{definition} Let $E \subseteq \mathbb{R}^d, d \geq 2$. The unoriented direction set $D(E)$ determined by $E$ is the image of the oriented direction set $\tilde{D}(E)$ under the double cover map 
$$\pi: S^{d-1} \to RP^{d-1}.$$ We may think of $D(E)$ as the set of parallel types of lines determined by $E$.
\end{definition}
 
Note for every set $E$, the oriented direction set $\tilde{D}(E) \subseteq S^{d-1}$ can be seen to be invariant under the antipodal map $a$. This is because given a direction achieved by a pair $x,y \in E$, the antipodal direction is achieved by the same pair under reversal of the roles of $x$ and $y$. Due to this it is easily seen that $\tilde{D}(E)=S^{d-1}$ if and only if $D(E)=RP^{d-1}$ and $\tilde{D}(E)$ is dense in $S^{d-1}$ if and only if $D(E)$ is dense in $RP^{d-1}$.

In this paper when we refer to the graph set of a function as a subset of $\mathbb{R}^d$ we will always mean a set, that up to rotation is the graph of a scalar-valued function $f: U \to \mathbb{R}$ where $U \subseteq \mathbb{R}^{d-1}$ is arbitrary with respect to the standard axis-system. The set $U$ is arbitrary and can even be empty or a single point. The function in this graph is said to be Lipschitz if there is a positive constant $C$ such that $|f(x)-f(y)| \leq C |x-y|$ for all $x,y \in U$, and non-Lipschitz otherwise. 

\vskip.125in 

Our main result is the following. 

\begin{theorem}
\label{theorem: main}
Let $E \subseteq \mathbb{R}^d$. Then there is a trichotomy in that exactly one of the following statements holds: \begin{itemize} 
\item i) $E$ is the graph of a Lipschitz function $f: A \to \mathbb{R}$ for some $A \subseteq \mathbb{R}^{d-1}$ (up to rotation) and $D(E)$ is not dense in $RP^{d-1}$. 

\item ii) $E$ is the graph of a non-Lipschitz function $f: A \to \mathbb{R}$ (up to rotation) and $D(E)$ is dense in $RP^{d-1}$ but is not equal to all of $RP^{d-1}$. 

\item iii) $E$ is not a graph of a scalar valued function and $D(E)=RP^{d-1}$. 
\end{itemize} 
\end{theorem}

\vskip.125in 

We can say a bit more with a few extra assumptions. 

\begin{corollary} 
\label{corollary: first} Let $E$ be a compact subspace of $\mathbb{R}^d$. Then there is the following trichotomy: \begin{itemize} 
\item i) $E$ is the graph of a Lipschitz function $f: A \to \mathbb{R}$ for some compact $A \subseteq \mathbb{R}^{d-1}$, up to rotation, and $D(E)$ is not dense in $RP^{d-1}$.
\item ii) $E$ is the graph of a continuous, non-Lipschitz function $f: A \to \mathbb{R}$, up to rotation, $A$ compact, and $D(E)$ is dense in $RP^{d-1}$ but is not equal to all of $RP^{d-1}$. 
\item iii) $E$ is not a graph of a scalar valued function and $D(E)=RP^{d-1}$. 
\end{itemize} 
\end{corollary}

Using the intermediate value theorem, one can further upgrade the result when $E$ is a compact, connected subset of $\mathbb{R}^2$. 

\begin{corollary} 
\label{corollary: second} 
Let $E$ be a compact, connected subset of $\mathbb{R}^2$. Then there is the following trichotomy: \begin{itemize} 
\item i) $E$ is the graph of a Lipschitz function $f: [a,b] \to \mathbb{R}$, up to rotation, and $D(E)$ is not dense in $RP^{1}$. 
\item ii) $E$ is the graph of a continuous, non-Lipschitz function $f: [a,b] \to \mathbb{R}$ (up to rotation) and $D(E)$ misses exactly one point of $RP^{1}$. 
\item iii) $E$ is not a graph of a scalar valued function and $D(E)=RP^{1}$. 
\end{itemize} 
\end{corollary}

Note the last corollary shows that once a compact, connected subset of the plane misses two directions, it misses a nonempty open set of directions and up to rotation, it is the graph of a Lipschitz map over a closed interval $[a,b]$ (where $a=b$ is a possibility). Also note, that it is known that when the vector space of continuous real-valued functions on the interval $[a,b]$, $a < b$ is given the sup-norm (uniform convergence norm), the nowhere differentiable continuous functions form a co-meagre set (topological analog of full measure subset). Thus ``most'' continuous functions are nowhere differentiable (and hence not Lipschitz as Lipschitz functions are almost everywhere differentiable), and hence ``most" continuous functions $f: [a,b] \to \mathbb{R}$ determine every possible secant slope in $\mathbb{R}$.

\vskip.125in 

We also establish the following:

\begin{proposition}
\label{proposition: countable}
Let $E \subseteq \mathbb{R}^d$ then $D(E)$ is countable if and only if $E$ is countable or $E$ is contained in a line, not necessarily through the origin. \end{proposition}
 
Putting these results together it follows that if $E \subseteq \mathbb{R}^{d}, d \geq 2$ and $E$ has Hausdorff dimension $> d-1$, then $D(E)$ is an uncountable dense subset of $\mathbb{R}P^{d-1}$. This is because the Hausdorff dimension of a graph set of a Lipschitz function $f: A \to \mathbb{R}, A \subseteq \mathbb{R}^{d-1}$ is at most $d-1$. It also follows immediately that if $E$ is a compact connected subset of the plane of Hausdorff dimension $>1$ then $D(E)$ misses at most one direction. 

\vskip.25in 
 
\section{Proofs of the main results} 

\vskip.125in 

Throughout the paper, $x$ and $y$ are $d$-dimensional vectors when we are in ${\Bbb R}^d$. In two dimensions, $(x,y)$ denotes a $2$-dimensional vector. 

\subsection{Proof of Proposition~\ref{proposition: countable}}

\vskip.125in 

When $E \subseteq \mathbb{R}^d$ is countable or contained in a line, it is clear that $D(E)$ is countable. So let us just prove the converse. Suppose $D(E)$ is countable and fix a point $x \in E$.
Then as $D(E)$ is countable, $E$ must be contained in a countable union of lines through $x$. If $E$ is countable we are done so assume it is uncountable. Then there must be a 
line $L$ through $x$ such that uncountably many elements of $E$ lie in this line. If there was another line $L'$ through $x$ that contained an element $e \in E$ besides $x$, we would obtain uncountably many directions generated by $E$ by noting that the lines through $e$ and the various uncountable elements of $E \cap L$ all have different directions. Thus it follows that if $E$ is uncountable, it must be contained in a single line.

\subsection{Proof of Theorem~\ref{theorem: main} }

\vskip.125in 

First consider $E \subset \mathbb{R}^d, d \geq 2$ with $D(E) \neq RP^{d-1}$. As $E$ fails to determine some direction, after a rotation, we may assume that $E$ fails to determine the 
$x_d$-axis direction. The projection $\pi: R^d \to R^{d-1}$ then determines a bijection $E \to \pi(E)=A \subseteq \mathbb{R}^{d-1}$. Let $f$ be the $d$-th coordinate function of the inverse 
of this bijection, it follows that $E=Graph(f)=\{ (x,f(x)) | x \in A \subseteq \mathbb{R}^{d-1} \}$ is the graph set of a scalar valued-function. Conversely when $E$ is the graph set of a 
scalar-valued function, by the vertical line test, it misses a direction.

Thus it follows that $D(E)=RP^{d-1}$ if and only if $E$ is not a graph set of a scalar-valued function. 

Thus it suffices for the remainder of the proof to only consider graph sets of functions 
$$f: A \to \mathbb{R} \ \text{with} \ A \subseteq \mathbb{R}^{d-1}.$$ 

Let us first consider the case $d=2$. Then $RP^1$, the space of lines through the origin can be identified as the circle $[0, \pi]/ (0 \sim \pi)$. This is because we can parameterize a line 
by the angle it makes with the $x$-axis, (angle $\pi$ and angle $0$ both give the same line, the $x$-axis). We can also use the slopes of these lines as the parametrization in which case $RP^1$ is identified with the one-point-compactification of $\mathbb{R}$ where the infinity slope corresponds to the $y$-axis. Under these identifications, dense subsets correspond to dense subsets so our conclusions are invariant no matter which picture we choose. 

Under the slope identification of $RP^1$, a graph set 
$$E=Graph(f)=\{ (x,f(x)): x \in A \subseteq \mathbb{R}^1 \}$$ has $D(E)$ given by the set of secant slopes 
$$ \left\{ \frac{f(y)-f(x)}{y-x}: x, y \in A, x \neq y \right\}$$ of the graph. 

Now suppose $E \subseteq \mathbb{R}^2$ does not have a dense set of directions, then $D(E)$ misses a nonempty open set of directions. We may rotate $E$ and assume this open set 
of directions includes the $y$-axis direction. Then as previously discussed it is the graph of a function $f: A \to \mathbb{R}$ whose secant slopes are bounded away from infinity as $D(E)$ 
avoids an open neighborhood of the infinite slope ($y$-axis).
 
In other words, there is a positive constant $C > 0$ such that 
$$\frac{|f(x)-f(y)|}{|x-y|} \leq C \ \text{for all} \ x, y \in A, x \neq y.$$ 

It follows immediately then that $f$ is Lipschitz and $E$ is the graph set of a Lipschitz function. Conversely, the graph of a Lipschitz function clearly misses an open set about the infinite slope direction via the same picture.

Thus we have established in the case of subsets $E \subseteq \mathbb{R}^2$ that $D(E)$ is not dense in $\mathbb{R}P^1$ exactly when $E$ is the graph set of a Lipschitz function 
and $D(E)=\mathbb{R}P^1$ exactly when $E$ is not the graph set of a scalar valued function. Hence Theorem~\ref{theorem: main} is proven in the case $d=2$.

Now we consider the general case when $d > 2$. If $D(E)$ is not dense, then after a rotation we can assume it misses an open set about the $x_d$-axis direction. In particular, there is an $\epsilon > 0$ such that it does not generate any line within angle $\epsilon$ of the $x_d$-axis. Again projection to the first $d-1$ coordinates is an injection on $E$ and so 
$E=Graph(f)$ where $f: A \to \mathbb{R}$, where $A \subseteq \mathbb{R}^{d-1}$ is the projection of $E$ to ${\Bbb R}^{d-1}$. 

Now fix an affine line $L$ (not necessarily through the origin) in $R^{d-1}$ then the unique affine $2$-plane in $\mathbb{R}^d$ which contains $L$ and a line parallel to the $x_d$-axis, intersects $E$ in the graph of the function $f$ restricted to $A \cap L$. This restriction of $f$ to the line $L$ is then Lipschitz with Lipschitz constant $\tan \left(\frac{\pi}{2}-\epsilon \right)$ by the 
argument in the $d=2$ case as its secant slopes are bounded away from the infinite slope which now corresponds to the $x_d$-axis direction. It is then clear that $f: A \to \mathbb{R}$ is Lipschitz as it is Lipschitz on $A \cap L$ for any line $L$, with a uniform Lipschitz constant $\tan(\frac{\pi}{2}-\epsilon)$. The Lipschitz condition involves two elements at a time and any two elements lie on a single affine line.

Conversely if $E=Graph(f)$ for some $f: A \to \mathbb{R}$, $A \subseteq \mathbb{R}^{d-1}$, then if $f$ is Lipschitz, there is a uniform bound on the secant slopes of 
$f$ restricted to $A \cap L$ for all lines $L$ through the origin in $\mathbb{R}^{d-1}$. This easily translates to the existance of an angle $\epsilon > 0$ such that $E$ does not generate any lines within $\epsilon$ angle of the $x_d$-axis.

Thus in the general $\mathbb{R}^d$ case we have seen that $D(E)=\mathbb{R}P^{d-1}$ if and only if $E$ is not the graph set of a scalar-valued function and 
$D(E)$ is not dense if and only if it is the graph set of a Lipschitz function. Thus the theorem is proven.

\subsection{Proof of Corollary~\ref{corollary: first} and Corollary~\ref{corollary: second}}

By the proof of Theorem~\ref{theorem: main}, when $D(E) \neq RP^{d-1}$, $E$ is (up to rotation) the graph set of a function $f: A \to \mathbb{R}$ where $A \subseteq \mathbb{R}^{d-1}$ 
is the projection of $E$ to $\mathbb{R}^{d-1}$. When $E$ is compact (connected), $A$ is also compact (connected). In particular when $d=2$, $A$ is a closed subinterval of $\mathbb{R}$ 
as every compact, connected subset of $\mathbb{R}$ is a closed interval (to see this just apply the intermediate and extreme value theorems to the inclusion map $i: A \to \mathbb{R}$.)

When this graph set $E$ is compact, it follows from the topological closed graph theorem (\cite{M75}) that $f$ is continuous. Finally when $d=2$, and $E$ is compact and connected, then (up to rotation) $E$ is the graph set of $f: [a,b] \to \mathbb{R}$. The secant formula 
$$\frac{f(x)-f(y)}{x-y}$$ defines a continuous function on the open triangle 
$$\{ (x,y) \in [a,b] \times [a,b]: x < y \}.$$ 

Since this triangle is connected it follows that  
$D(E)$ is connected in $RP^{1}$. When the function $f: [a,b] \to \mathbb{R}$ is continuous but not Lipschitz, it follows from Theorem~\ref{theorem: main} that its secants are dense and not bounded, this combined with $D(E)$ is connected is enough to conclude then that 
$$D(E)=RP^{1}-\{ \text{ y-axis direction } \}$$ using the intermediate value theorem. Use the identification of $RP^1$ via slopes so that 
$$RP^1-\{ \text{ y-axis direction } \}$$ is identified with $\mathbb{R}$. Thus all directions are achieved except the $y$-axis direction in this case. This concludes the proofs of these corollaries.

\end{document}